\documentclass[twoside,11pt]{article}
\usepackage{amssymb,amsmath,amsthm,stmaryrd, amsopn}
\usepackage{amscd} 
\usepackage{psfig,latexsym,graphicx}
\usepackage{fullpage}
\newenvironment{Quote}
{\begin{list}{}{%
   \vspace{-.8cm}
   \setlength{\rightmargin}{7mm}
   \setlength{\leftmargin}{7mm}}
    \item[]}
{\end{list}}

\refstepcounter{equation}

\newtheorem{lem}{LEMMA}
\newtheorem{theo}[lem]{THEOREM}

\newtheorem{hyp}[lem]{\indent Hypothesis}

\theoremstyle{definition}

\newtheorem{rem}[lem]{\indent Remark}

\newtheorem{ex}[lem]{\indent Example}
\refstepcounter{lem}

\renewcommand{\descriptionlabel}[1]%
     {\hspace{\labelsep}\textsf{#1}}

\newtheorem{alem}[lem]{LEMMA\, A\hspace{-.15cm}}

\newtheorem{blem}[lem]{LEMMA\, B\hspace{-.15cm}}

\newtheorem{bprop}[lem]{PROPOSITION\, B\hspace{-.15cm}}
\newtheorem{brem}[lem]{Remark\, B\hspace{-.15cm}}

\newcommand{\mycaption}[1]{\begin{quote}
 \caption{\small#1\normalsize}\end{quote}\vspace{-3em}}

\newcommand\Z{{\mathbb Z}}
\newcommand\R{{\mathbb R}}

\newcommand\A{{\mathcal A}}
\newcommand\B{{\mathcal B}}

\newcommand\cC{{\mathcal C}}

\newcommand\ssm{{\smallsetminus}}
\font\bit=cmssi12 at 11truept

\def\Ly{{\rm Lyap}}

\begin{document}

\markboth{\centerline{\sc Schwarzian Derivatives and Cylinder Maps}}
{\centerline{\sc Araceli Bonifant and  John Milnor}}

\title{\bf Schwarzian Derivatives and Cylinder Maps}

\author{\bf {Araceli Bonifant
\,\, and \,\, John Milnor
}}

\setcounter{footnote}{1}
\setcounter{equation}{0}

\date{}

\maketitle

\thispagestyle{empty}
\def\IMSmarkvadjust{0 pt}
\def\IMSmarkhadjust{0 pt}
\def\IMSmarkhpadding{0 pt}
\def\IMSpubltext{Published in modified form:}
\def\SBIMSMark#1#2#3{
 \font\SBF=cmss10 at 10 true pt
 \font\SBI=cmssi10 at 10 true pt
 \setbox0=\hbox{\SBF \hbox to \IMSmarkhpadding{\relax}
                Stony Brook IMS Preprint \##1}
 \setbox2=\hbox to \wd0{\hfil \SBI #2}
 \setbox4=\hbox to \wd0{\hfil \SBI #3}
 \setbox6=\hbox to \wd0{\hss
             \vbox{\hsize=\wd0 \parskip=0pt \baselineskip=10 true pt
                   \copy0 \break%
                   \copy2 \break%
                   \copy4 \break}}
 \dimen0=\ht6   \advance\dimen0 by \vsize \advance\dimen0 by 8 true pt
                \advance\dimen0 by -\pagetotal
	        \advance\dimen0 by \IMSmarkvadjust
 \dimen2=\hsize \advance\dimen2 by .25 true in
	        \advance\dimen2 by \IMSmarkhadjust

%
%
  \openin2=publishd.tex
  \ifeof2\setbox0=\hbox to 0pt{}
  \else 
     \setbox0=\hbox to 3.1 true in{
                \vbox to \ht6{\hsize=3 true in \parskip=0pt  \noindent  
                {\SBI \IMSpubltext}\hfil\break
                \input publishd.tex 
                \vfill}}
  \fi
  \closein2
  \ht0=0pt \dp0=0pt
 \ht6=0pt \dp6=0pt
 \setbox8=\vbox to \dimen0{\vfill \hbox to \dimen2{\copy0 \hss \copy6}}
 \ht8=0pt \dp8=0pt \wd8=0pt
 \copy8
 \message{*** Stony Brook IMS Preprint #1, #2. #3 ***}
}

\SBIMSMark{2006/7}{October 2006}{}

\begin{abstract}
We describe the way in which the sign of the Schwarzian 
derivative for a family of diffeomorphisms  of the interval $I$ affects the 
dynamics of an associated many-to-one
skew product map of the cylinder $(\R/\Z)\times I$.
\end{abstract}

\vspace{.2cm}

\vspace{.2cm}

\noindent
\setcounter{lem}{0}

\section{Introduction.}
Ittai  Kan has described a simple example of a skew product
 map from the cylinder
$(\R/\Z)\times I$ to itself such that the two boundary circles are
 measure theoretic attractors whose attracting
basins are {\bit intermingled\/}, in the sense that the intersection of
any nonempty open set with either basin has strictly positive measure.
(See  \cite{Kan}.)
This note will consist of three variations on the maps which he introduced.

Sections \ref{NS} and 4 will describe Kan's example in slightly more
generality, emphasizing the importance of
{\bit negative Schwarzian derivative\/}. Section \ref{SP} will
show that if we substitute {\bit positive Schwarzian derivative\/} then the 
behavior will change drastically, and almost all orbits will have a common 
asymptotic distribution. In the case of {\bit zero Schwarzian derivative,}~
\S\ref{ZS}  will prove in some cases (and conjecture in others) that
almost all orbits spend most of the time extremely close to one of the two
 cylinder boundaries; but that each such orbit passes
 from the $\epsilon$-neighborhood of one boundary circle
to the $\epsilon$-neighborhood of the other infinitely many times
on such an irregular schedule that there is no asymptotic measure.

Most technical details  are relegated to the two appendices.

\begin{figure}
\centerline{\psfig{figure=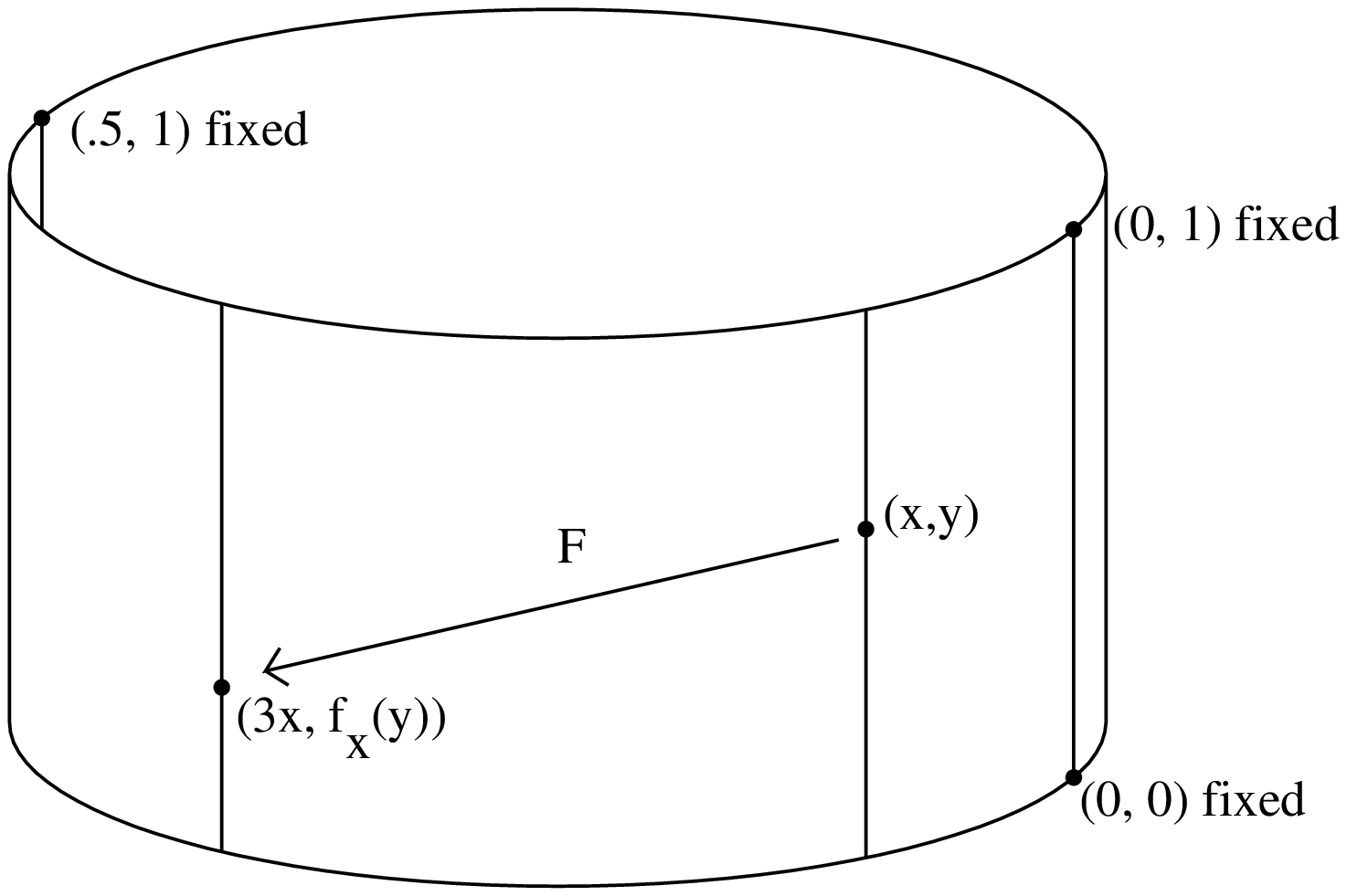,height=1.8in}}
\mycaption{\bit  The cylinder map $F$ in the case $~k=3$.}
\end{figure}

\noindent
\setcounter{lem}{0}
\section{\bf Preliminaries.}\label{Preli} Let $I=[0,1]$, and let
$\cC$ be the cylinder $(\R/\Z)\times I$ with boundaries 
$\A_0=(\R/\Z)\times\{0\}$ and \break $\A_1=~(\R/\Z)\times\{1\}$.
Let $~F:\cC\to \cC~$ be a $~C^3$-differentiable map of the form
\begin{equation}\label{eq-for}
 F(x,y)=\big(kx,\,f_x(y)\big)\,, 
\end{equation}
where $k\ge 2$ is a fixed integer, and where each $f_x:I\to I$ is a
diffeomorphism with $f_x(0)=0$ and  $f_x(1)=1$. Thus the derivative
$$f_x'(y)=\partial f_x(y)/\partial y$$
must be strictly positive everywhere. 
\smallskip

We next introduce two key concepts that will be needed.

\begin{quote}
{\bf LEMMA A.1.}
{\it For $\iota$ equal to zero or one, let $\B_{\iota}$ be the
attracting basin of the circle $\A_\iota$. If the transverse Lyapunov exponent
\begin{equation}
\Ly(\A_{\iota})=\int\limits_{\R/\Z} \log\left(f'_x(\iota)\right)dx 
\end{equation}
is negative, then the basin $\B_{\iota}$ has strictly positive measure.
In fact, for almost every $x\in \R/\Z$ the basin $\B_\iota$ intersects
the ``fiber'' $x\times I$ in an interval of positive length. On the other
hand, if $~\Ly(\A_{\iota})>0~$ then $\B_{\iota}$ has measure zero.}
\end{quote}

\noindent The proof will be given in Appendix A. 
\medskip 

In fact, whenever $~\Ly(\A_\iota)<0~$ it is not hard to see that the circle 
$~\A_{\iota}~$ is a {\bit measure attractor\/}. By this we mean
 that  it satisfies the following two conditions:
\begin{enumerate}
\item[{\bf 1}.] $\A_\iota$  is a {\bit minimal measure attracting set\/}, that
is, it has an attracting basin of positive measure,
 but no closed proper subset has a basin of positive measure.
\item[{\bf 2.}] Furthermore,
 $\,\A_\iota\,$ contains a dense orbit, and hence cannot be expressed
 as the union of strictly smaller closed invariant sets.\footnote{The following
example shows that Condition {\bf 2} does not follow from Condition {\bf 1}.
Consider a flow in the plane such that all orbits near infinity spiral in
towards a
figure eight-curve, while all orbits inside either lobe of the figure-eight
spiral in towards an attracting equilibrium point. Then the figure-eight is
a minimal measure attracting set with no dense orbit.}
\end{enumerate}
\medskip

Recall that the {\bit Schwarzian derivative\/} of an interval 
$C^3$-diffeomorphism $f$ is defined by the formula
\begin{equation}\label{eq-ns}
{\mathcal S} f(y)~=~ \frac{f'''(y)}{f'(y)}  - \frac{3}{ 2}\left(\frac{f''(y)}{ 
 f'(y)}\right)^2.
\end{equation}

\noindent
We will make a particular study of
maps $~F(x,y)=\big(kx,\,f_x(y)\big)~$ such that the 
Schwarzian $~{\mathcal S} f_x(y)~$ has constant sign for almost all 
$~(x,y)\in\cC$. Maps $~f_x~$ with $~{\mathcal S}(f_x) <0~$ almost everywhere
have the basic property of
{\it increasing\/}  the cross-ratio $~\rho(y_0,y_1,y_2,y_3)~$  for all
$~y_0 < y_1 < y_2 < y_3~$ in the interval. (See Appendix B.) Similarly, if 
$~{\mathcal S}(f_x) >0~$ (or if  $~{\mathcal S}(f_x)\equiv 0$),
then $~f_x~$ will decrease (or will preserve) all such cross-ratios.

\medskip

\setcounter{lem}{0}

\section{\bf  Negative Schwarzian.}\label{NS}

\begin{Quote}
\begin{lem}\label{lem-ns}
If $~{\mathcal S} f_x(y)~$ has constant sign $($positive, negative or, zero$)$
for almost all $~(x,\,y)$,~  then
$~\Ly(\A_0) + \Ly(\A_1)~$ has the same sign. In particular, if
$~{\mathcal S} f_x(y)<0~$ for almost all $~(x,\,y)$,~  then
\begin{equation}\label{eq-sumlya}
\Ly(\A_0) + \Ly(\A_1) < 0\,,
\end{equation}
hence at least one of the two boundaries has a basin of positive measure.
\end{lem}
\end{Quote}
\smallskip

{\bf Proof.\/} Lemma B\ref{lem-prod} (in Appendix B) will show that
 $f'_x(0)f'_x(1) < 1$ whenever $~f_x~$ has negative Schwarzian.
Integrating the logarithm of this inequality over $\R/\Z$, the inequality
  (\ref{eq-sumlya}) follows. Thus
the transverse Lyapunov exponent is negative for at least one of the two
boundaries. Hence the associated basin has positive measure
by Lemma A\ref{A.1.}.\qed
\medskip

\begin{Quote}
\begin{theo}\label{th-fm}
 If $\,{\mathcal S} f_x(y)<0\,$ almost everywhere, and if both basins
have positive measure,\footnote{We don't know whether this hypothesis is
necessary.}
then  $\B_0\cup \B_1$ has full measure. In fact,
there is an almost everywhere defined measurable function 
$\;\sigma:\R/\Z\to I\;$ such that
\begin{eqnarray*}
 (x,y)\in \B_0\qquad &{\it whenever}&\qquad y<\sigma(x)\,,\cr
 (x,y)\in \B_1\qquad &{\it whenever}&\qquad y>\sigma(x)\,.\cr
\end{eqnarray*}\vskip -.5cm
More generally, the same statement is true if the $k$-tupling map on the circle
is replaced by any continuous ergodic transformation $~g~$
on a compact space with $~g$-invariant probability measure.
\end{theo}
\end{Quote}
\smallskip

In fact we will usually consider maps $f_x(y)$ for which the
behavior of $F$ near the two boundaries is similar enough so that
$~\Ly(\A_0)~$ and $~\Ly(\A_1)~$ are equal to each other (or at least
 have the same sign). For such maps, the condition
$~{\mathcal S} f_x<0~$ will guarantee that both attracting basins
have positive measure.
\medskip

{\bf Proof of Theorem \ref{th-fm}.\/} Since each $f_x$ is an orientation preserving homeomorphism,
there are unique numbers $0 \leq \sigma_0(x) \leq 
\sigma_1(x) \leq 1~$ defined by the  property that the orbit of $(x,y)$:
 \begin{center}
converges to $\A_0$\quad if \quad $y\, < \,\sigma_0(x)$\\
converges to $\A_1$\quad if \quad $y\, > \, \sigma_1(x)$\\
\hspace{-1.4cm}does not converge to either circle\quad  if\quad
$\sigma_0(x)\, <\, y\, <\, \sigma_1(x)\,.$
\end{center}
\noindent 
Thus, the area of $\B_0$ can be defined as $\int\sigma_0(x)\,dx$.
Since this is assumed to be positive, it follows that the set of all
$x\in\R/\Z$ with $\sigma_0(x)>0$ must have
positive measure. In fact, the evident identity
$~\sigma_0(kx)=f_x\big(\sigma_0(x)\big)~$ implies that
this set is fully invariant under the ergodic map $x\mapsto kx$.~ Hence it
 must actually have full measure. Similarly, the
set of $x$ with $\sigma_1(x)<1$ must have full measure.\smallskip

We will make use of the property that a map $f_x$
 of negative Schwarzian derivative increases the cross-ratio
$$\rho(0,\,y_1,\,y_2,\,1)~=~\frac{y_2\,(1-y_1)}{y_1\,(1-y_2)}\,,$$
that is:
$$\rho\big(0,\,f_x(y_1),\, f_x(y_2),\,1\big)~>~
 \rho(0,\,y_1,\,y_2,\,1)~>~1
\qquad{\rm for~~~all}\qquad 0<y_1<y_2<1\,.$$
(See Lemma B\ref{lem-nsi}.) Suppose that the inequalities
 $0<\sigma_0(x)<\sigma_1(x)<1$ were true for a set of $x\in\R/\Z$ of
positive Lebesgue measure, then the function
$$ r(x)~~=~~ 
\rho\big(0,\,\sigma_0(x),\,\sigma_1(x),\,1\big)~~\ge~~1$$
would satisfy  $r(kx)>r(x)$ on a set of positive measure, with $r(kx)\ge r(x)$
everywhere. It would follow that
$$\int\limits_{\R/\Z}\frac{dx}{r(kx)}~~<~~\int\limits_{\R/\Z}\frac{dx}{r(x)}
\,.$$
But this is impossible: Lebesgue measure is invariant under push-forward
by the map $x\mapsto kx$, and it follows that
 $\int\phi(kx)\,dx=\int\phi(x)\,dx$ for any bounded measurable function $\phi$.
This contradiction proves that we must have
$\sigma_0(x)=\sigma_1(x)$ almost everywhere; and we define $\sigma(x)$ as this
common value.
\qed
\medskip

\begin{rem}\label{rem-sepm}
We can then define the {\bit separating measure\/} $\beta$ on $\cC$ to be the
 push-forward, under the section,  $x\mapsto (x,\,\sigma(x))$, of the Lebesgue
measure $\lambda_x$, on $\R/\Z$. Evidently  $\beta$ is an $F$-invariant 
ergodic probability measure which in some sense
describes the ``boundary'' between the two basins. Since $~0<\sigma(x)<1~$
almost everywhere, it follows easily that both boundaries have measure
$~\beta(\A_\iota)=0$.
\end{rem}

\setcounter{lem}{0}

\section{\bf Intermingled Basins.} Now assume the  following.

\begin{Quote}
\begin{hyp}\label{hyp-ib}
There exist angles $x^-$ and $x^+$ in $\R/\Z$, both
fixed under multiplication by $k$,
 such that $f_x(y)< y$ for all $0<y<1$ and all $x$ in a neighborhood of
$x^-$, and such that $f_x(y)>y$ for all $0<y<1$ and all $x$ near $x^+$.
\end{hyp}
\end{Quote}

It follows that
the entire vertical line segment $\{x^-\}\times[0,1)$ is contained in
the basin $\B_0$, and  that the entire segment
 $\{x^+\}\times(0,1]$ is contained in the basin $\B_1$.
\medskip

\begin{Quote}
\begin{theo}\label{th-inm}{\bf (Intermingled Basins).\/} 
If Hypothesis \ref{hyp-ib}
is satisfied, and if both basins have positive measure,
then the two basins are {\bit intermingled\/}. That is, for every
non-empty open set $U\subset \cC$,
both intersections $\B_0\cap U$ and $\B_1\cap U$ have strictly positive 
measure.
\end{theo}
\end{Quote}
\smallskip

{\bf Proof.\/} Define measures $\mu_0$ and $\mu_1$ on the cylinder by setting
$\mu_\iota(S)$ equal to the Lebesgue measure of the intersection
$\B_\iota\cap S$ for $\iota$ equal to zero or one and for any measurable
set $S$.  Clearly the {\bit support\/} ${\bf supp}(\mu_\iota)$,
that is the smallest closed set which has full measure under $\mu_\iota$,
is fully $F$-invariant. We must prove that this support
is equal to the entire cylinder.

To begin, choose any point $(x_0, \,y_0)\in{\bf supp}(\mu_0)$ with $0<y_0<1$.
Construct a backward orbit
$$ \cdots~\mapsto~ (x_{-2},\,y_{-2})~\mapsto~ (x_{-1},\,y_{-1})~\mapsto~
 (x_{0},\,y_{0})$$
under $F$ by induction,
letting each $x_{-(k+1)}$ be that preimage of $x_{-k}$ which is closest to
$x^-$. Then it is not difficult to see that this backwards sequence
converges to the point $(x^-,\,1)$. Since ${\bf supp}(\mu_0)$ is closed
and $F$-invariant, it follows that $(x^-,\,1)\in{\bf supp}(\mu_0)$. But
the iterated pre-images of $(x^-,\,1)$ are everywhere dense in the
upper boundary circle $\A_1$, so $\A_1$ is contained in ${\bf supp}(\mu_0)$.
Since the basin $\B_0$ is a union of vertical line segments $x\times
 \big[0,\,\sigma_0(x)\big)$ or $x\times\big[0,\,\sigma_0(x)\big]$, it follows easily that
${\bf supp}(\mu_0)$ is the entire cylinder.

 The proof for $\mu_1$ is completely analogous.\qed
\medskip

\begin{rem}
 In place of a fixed point on the circle, we could equally well
use a periodic point $k^px\equiv x$ $({\rm mod}\, \Z)$. It is only necessary to
check that the iterated map $F^{\circ p}$ satisfies the required hypothesis.
\end{rem} 
\smallskip

\begin{figure}[h]
\centerline{\psfig{figure=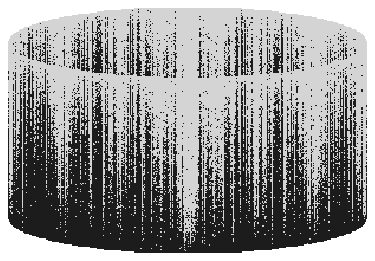,height=1.8in}}
\mycaption{\bit  Intermingled basins for the cylinder map $F$ of Example 
\ref{ex-kan}.\label{f-kan}}
\end{figure}
\smallskip

\begin{ex}\label{ex-kan}
 Following  \cite{Kan}, let
\begin{equation}\label{eq-kan}
q_a(y)=y\,+\,ay(1-y)\,.
\end{equation}
If $~|a|< 1$,~
then $~q_a~$ maps the unit interval diffeomorphically onto itself, with
$q_a(0)=0$ and $q_a(1)=1$. It is easy to check that $\,{\mathcal S} q_a(y) 
< 0\,$
whenever $a\ne 0$. It then follows from Lemma B\ref{lem-nsi} that
$\,q_a(y)$ has the property of  increasing cross-ratios.

Choose $0<\epsilon<1$, and let $~p(x)=\epsilon\,\cos(2\pi x)$.~
Then for any $k\ge 3$ the map
$$ F(x,y)~=~\big(kx,\,f_x(y)\big) \qquad{\rm where}\qquad
 f_x(y)=q_{p(x)}(y)$$
will satisfy Hypothesis \ref{hyp-ib} and also the  hypotheses of 
Theorem~\ref{th-fm}. In fact, we can take 
$x^+=0$, and choose $x^-$ to be a fixed point which lies 
between $1/3$ and $2/3$.  For example, take
$$x^- =
\begin{cases} 
1/2\,, & \quad {\rm for}\quad k~~~ {\rm odd},\cr
k/(2k-2)\,, & \quad {\rm for}\quad k\ge 4~~~ {\rm even}.
\end{cases}
$$
{\it Thus we obtain explicit examples of maps with intermingled basins.\/}
(Compare Fig.~\ref{f-kan}.) \smallskip
 
(In fact this argument will work for $~k=2~$ also, using the periodic orbit
$1/3\leftrightarrow 2/3$ in place of a fixed point.)
\end{ex}
\smallskip

\begin{rem}
Very similar examples of intermingled basins can be observed in rational
maps of the projective plane. (Compare \cite[\S6]{BDM}.)
It would be very interesting to know to what extent the examples in the
following two sections also have analogs among such rational maps.
\end{rem}

\setcounter{lem}{0}
\section{Positive Schwarzian}\label{SP}
In this section we continue to study the cylinder maps
$~F(x,y)=(kx,\, f_x(y))\,$,~ but now assume that $~{\mathcal S} f_x>0~$ almost 
everywhere.
\smallskip

\begin{Quote}
\begin{theo}\label{theo-meas}
If $\,{\mathcal S} f_x(y)>0\,$ for almost all
$(x,y)$, then at least one of the transverse Lyapunov exponents
$\Ly(\A_0)$ and $\Ly(\A_1)$ is strictly positive. If both are strictly 
positive, then $F$ has an asymptotic measure.\footnote{Terms such as: natural 
measure or physical measure  are also used in the literature to 
denote this type of measure.} That is, there is a uniquely defined probability
 measure $~\nu~$ on the cylinder such that, for Lebesgue almost every orbit  
 $~(x_0,y_0)\mapsto (x_1,y_1)\mapsto \cdots~$ and for every continuous test 
function $\chi:\cC\to\R$, the time averages
$$ \frac{1}{n}\Big(\sum_{i=0}^{n-1} \chi(x_i,\,y_i)\Big) $$
converge to the space average $~\int_\cC \chi(x,y)\,d\nu(x,y)~$ as
$~n\to\infty$.~ $($Briefly, almost every orbit is
{\bit uniformly distributed\/} with respect to $\nu.)$
Furthermore, both boundaries of $\cC$
have asymptotic measure $~\nu(\A_0)=\nu(\A_1)~$ equal to zero.
\end{theo}
\end{Quote}
\smallskip

\noindent Thus, under these hypotheses, almost all orbits of $~F~$ have the
same asymptotic distribution.
\medskip

{\bf Outline of the Proof.\/} Since the proof of this theorem will be slightly 
circuitous,  we first outline the main steps.

\begin{itemize}
\item[$\bullet$] First the circle $~{\mathbb R}/{\mathbb Z}~$ of the 
previous section will be replaced by the solenoid
$$\Sigma~=~\Sigma_k~=~\lim_{\leftarrow}\, ({\mathbb R}/k^n{\mathbb Z})~,$$
and the many-to-one map $~F~$ of $~({\mathbb R}/{\mathbb Z})\times I~$ will 
be replaced by the associated invertible map $~\widetilde{F}~$ from $~\Sigma \,
\times\, I~$ to itself.

\item[$\bullet$] Since $~\widetilde{F}^{-1}~$ has negative Schwarzian on 
each fiber, Theorem~\ref{th-fm} asserts that the union of the attracting
basins of the two boundaries $~\Sigma \times \{0\}~$ and $~\Sigma \times \{1\}~$ under 
the map $~\widetilde{F}^{-1}$ will have full measure. In fact there is an
almost everywhere defined section
\begin{equation}\label{e:s1}
\tilde{x}\, \mapsto\, (\,\tilde{x},\, \sigma(\tilde{x})\,)
\end{equation}
from $~\Sigma~$ to $~\Sigma\, \times \,I~$ which ``separates''
 the two attracting basins.

\item[$\bullet$] There is a standard ergodic invariant probability 
measure $~\mu_{\Sigma}~$ on the solenoid. Pushing it up to the graph of
 $~\sigma~$ under
the section (\ref{e:s1}), we obtain an ergodic invariant 
probability measure $~\widetilde{\nu}~$ on $\Sigma\, \times\, I$.
\item[$\bullet$] Since almost all points are pushed {\it away} from the graph
of $\sigma$ by the inverse map $\widetilde{F}^{-1}$, it follows that they are
pushed {\it towards} this graph by the map $\widetilde{F}$. In this way, we
see that $~\widetilde{\nu}~$ is an asymptotic measure for $~\widetilde{F}~$.
\end{itemize}

\begin{itemize}
\item[$\bullet$] Finally, we denote by  $~\nu~$ the push-forward of 
$~\widetilde{\nu}~$ under the projection
$$\Sigma\, \times \,I\, \rightarrow \,({\mathbb R}/{\mathbb Z})\, \times\, I.$$
This will be the required asymptotic measure
 for the original cylinder map $~F$.
\end{itemize} 
\smallskip

\begin{figure}[h]
\centerline{\psfig{figure=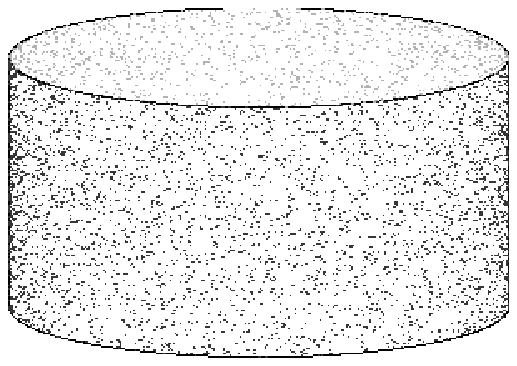,height=1.8in}}
\mycaption{\bit  50000 points of a randomly chosen
 orbit for the cylinder map $F$ of Example 
\ref{ex-invkan}.\label{f:invq-rand}}
\end{figure}
\smallskip

\begin{ex}\label{ex-invkan}
One example of a family of interval diffeomorphisms with positive
Schwarzian is given by the inverses
$$f_x(y)~=~q_{p(x)}^{-1}(y)\,,$$
where $q_{p(x)}(y)$ is the quadratic map (\ref{eq-kan})
of Example~\ref{ex-kan}. Here $\,{\mathcal S} f_x(y) > 0\,$ whenever
 $p(x) \neq 0$.
 (Compare Proposition B\ref{prop-propt} in the appendix.)
For this special example, the asymptotic measure $~\nu~$ turns out to be
 precisely equal to the
standard Lebesgue measure $~\lambda~$ on the cylinder.  In other words:

\centerline{\it Randomly chosen orbits are uniformly distributed with respect 
to Lebesgue measure.\/}\smallskip

\noindent
(Compare Fig.~\ref{f:invq-rand}). To prove this statement, one
only needs to show that {\bit Lebesgue measure is 
$~F$-invariant}.  In fact
there are $~k~$ branches of $~F^{-1}~$ on any small open set $~U \subset \cC$,
each given by
\begin{equation}\label{F^-1}
F^{-1}(x,y)~=~\big(~x/k,~y+ \epsilon\,\cos(2\pi x/k)y(1-y)\big)
\end{equation}
for one of the $~k~$ choices of $~x/k~~({\rm mod}~{\mathbb Z})$.~
The Jacobian of this branch (\ref{F^-1}) is equal to
$$~\big(1+\epsilon\,\cos(2\pi x/k)(1-2y)\big)/k\,.$$
Since the sum of
$\cos(2\pi x/k)~=~\Re e^{2\pi i x/k} $
over the $~k~$ choices for $~x/k~$ is zero, the sum of Jacobians is $+1$,
which means that $~F~$ preserves the Lebesgue measure $~\lambda$.
Now Theorem \ref{theo-meas} asserts that an asymptotic measure exists. Such a
 measure $~\nu~$ must necessarily be equal to the weak limit of
 $~\left(\lambda+F_*\lambda+\cdots+ F_*^{\circ n-1}\lambda\right)/n~$
as $~n\to\infty$;~ and it follows that
 $~\nu~$ is precisely equal to $~\lambda$.
\end{ex}
\smallskip

{\bf Proof of Theorem \ref{theo-meas}.} The argument begins as follows.
Denote by $\Sigma$ the solenoid of backwards orbits
$$\widetilde{x}=\{\cdots\mapsto x_{-2} \mapsto x_{-1} \mapsto x_0 \}=
\{x_{-n}\}$$
under the map $x \mapsto kx$. Thus $\Sigma$ maps homeomorphically onto itself
under multiplication by $k$, with the right
shift map $$\{\cdots \mapsto x_{-2}\mapsto x_{-1}\mapsto x_0\} \mapsto
\{\cdots\mapsto x_{-3}\mapsto x_{-2}\mapsto x_{-1}\}$$
as inverse. There is a standard probability measure $\;\mu_{\Sigma}\;$
on $\Sigma$, defined by the requirement
that each projection $~\widetilde x\mapsto x_{-n}~$ is measure preserving.
\smallskip


\begin{Quote}
\begin{lem}
$~~{\|\widetilde{x}\|/k} ~\leq~ \| k\widetilde{x}\|
~\leq~ k\, \| \widetilde{x}\| ~~$
for all $~\widetilde{x} \in \Sigma.$
\end{lem}
\end{Quote}
\smallskip

{\bf Proof.} This follows easily from the definition.\qed
\medskip

Given maps $f_x(y)$ as in \S \ref{Preli}, 
we again consider the associated map $~F:\cC\to\cC$.
If $\,{\mathcal S} f_x(y)>0\,$ almost everywhere, then by the analogue to
Lemma \ref{lem-ns}, $$\;\Ly(\A_0)+\Ly(\A_1)>0.$$
There is also a {\bit natural extension\/}
 $~\widetilde{F}:\,\Sigma\,\times\, I~\to~\Sigma\,\times\,I$
of the map $F$. This is a homeomorphism defined by the formula
\begin{equation}\label{eq-calF}
 \widetilde{F}(\widetilde x,\;y)~=~\big(k\,\widetilde x,\; f_{x_0}(y)\big)\,.
\end{equation}
\smallskip

\begin{Quote}
\begin{lem}\label{lem-gam}
If $\,{\mathcal S} f_x(y) > 0\,$ for almost all $(x,y)$, and if both 
$\Ly(\A_0)$ and $\Ly(\A_1)$  are  strictly positive, then there exists a 
measurable  function $\sigma: \Sigma \rightarrow I$, defined almost everywhere,
 and satisfying the identity
$$\sigma(k\widetilde{x})=f_{x_0}\big(\sigma(\widetilde{x})\big)$$
for almost all $~\widetilde{x} \in \Sigma$. It follows that
the {\bit graph} of $\sigma$, that is the set of all pairs
$$(\widetilde{x},\, \sigma(\widetilde{x}))\; \in\; \Sigma \;\times\; I,$$
is invariant under the extended map $~\widetilde{F}:\Sigma\times I\to\Sigma
\times I$, so that $~\widetilde{F}\big(\widetilde{x},\,
\sigma(\widetilde{x})\big)=\big(k\widetilde{x},\, \sigma(k\widetilde{x})\big)~$
for almost all $\widetilde x$.
\end{lem}
\end{Quote}\smallskip

{\bf
Caution:} In cases of interest, this function $~\sigma~$ will not be continuous
and will not be everywhere defined.
\medskip

{\bf Proof of Lemma \ref{lem-gam}.\/} 
We apply Theorem \ref{th-fm} to the inverse map
$~\widetilde F^{-1}$,~ with the $~k$-tupling map on the circle replaced by the
right shift map on the solenoid. This yields a measurable section
$~\widetilde x~\mapsto~\big(\widetilde x,\, \sigma(\widetilde x)\big)~$
from $~\Sigma~$ to $~\Sigma\times I$.~ The required $~F$-invariance property
then follows easily.
\qed\medskip

Next we will show that almost every orbit under $~\widetilde{F}$ converges,
in a suitable sense, to the graph of $~\sigma$. Recall that
$~\sigma(\widetilde x)~$ is well defined~ and belongs to the open interval
$~(0,\,1)$~ for almost every $~\widetilde x\in\Sigma$. Thus, for
almost every point $~(\widetilde x,\, y)\in\Sigma\times (0,1)$,~ the quantity
$$  r(\sigma(\widetilde x),\,y)~=
~\big|\log\rho\big(0,\,\sigma(\widetilde x),\,y,\,1\big)\big|~\ge~0  $$
is defined and finite, vanishing if and only if $~y=\sigma(\widetilde x)$.
We will think of $~r(\sigma(\widetilde x),\,y)~$ as a measure of distance 
between $~\sigma(\widetilde x)~$ and $~y$.
\medskip

\begin{Quote}
\begin{lem}\label{lem-rto0}
Under the same hypothesis, for almost every orbit
$~(\widetilde x_0,\,y_0)\mapsto(\widetilde x_1\,y_1)\mapsto\cdots~$ under 
$~\widetilde{F}~$ this measure of distance
 $~r(\sigma(\widetilde x_n),\,y_n)~$ converges to zero as $~n\to\infty$.
\end{lem}
\end{Quote}
\smallskip

{\bf Proof.} Since the map $~\widetilde{F}(\widetilde x,\,y)=
\big(k\widetilde x,\, f_{x_0}(y)\big)~$  decreases cross-ratios on each fiber
(compare Lemma B\ref{lem-nsi}), we have
$$ r\big(\sigma(k\widetilde x),\,f_{x_0}(y)\big)~<~r\big(\sigma(\widetilde x),
\,y\big)$$
almost everywhere. For any constant $r_0>0$, let $N(r_0)$ be the strip
consisting of all $$~(\widetilde x,\,y)\in\Sigma\times(0,1)\qquad{\rm 
with} \qquad r(\sigma(\widetilde x),\,y)~\le~r_0\,.$$
Evidently $N(r_0)$ is mapped into itself by $~\widetilde{F}$.
Given constants  $~0<r_0<r_1$,~ we will also consider the difference set
$N(r_1)\ssm N(r_0)$. Let  $~s(\widetilde x)\le 1~$ 
be the supremum of the ratio
\begin{equation}\label{rat}
\frac{r\big(\sigma(k\widetilde x),\,f_{x_0}(y)\big)}{r(\sigma(\widetilde x),
\,y)}\qquad {\rm for}\qquad (\widetilde x,\,y)\in N(r_1)\ssm N(r_0)\,,
\end{equation}
or in other words for 
$~r_0< r(\sigma(\widetilde x),\,y)\le r_1$. Since the Schwarzian is positive 
almost everywhere, it is not hard to see that this supremum
satisfies $~s(\widetilde x)<1~$ for almost all $\widetilde x$.
(Here we make essential use of the fact that $~r_0>0$,~ since if
 the Schwarzian vanishes at $(\widetilde x,\,\sigma(\widetilde x))$ then the
ratio (\ref{rat}) would tend to $1$ as $y$ tends to $\sigma(\widetilde x)$.)

Therefore the average of $~\log s(\widetilde x)~$ over the solenoid
is strictly negative.
A straightforward application of the Birkhoff Ergodic Theorem then shows
that, for almost every $~\widetilde x_0\mapsto\widetilde x_1\mapsto\cdots~$, 
some partial product of the $~s(\widetilde x_j)~$ satisfies
$$ s(\widetilde x_0)\cdots s(\widetilde x_{n-1})~<~r_0/r_1\,.$$
This means that the iterate $~\widetilde{F}^{\circ n}~$ maps $N(r_1)$ into 
$N(r_0)$. Since $~0<r_0<r_1~$ can be arbitrary, this
completes the proof of Lemma \ref{lem-rto0}.\qed
\bigskip

Now define the probability measure $~\widetilde{\nu}~$ on $~\Sigma\,\times\,
 I~$ to be the push-forward of the standard measure $~\mu_{\Sigma}~$
on the solenoid under this section $~\widehat\sigma:\widetilde x\mapsto
\big(\widetilde x,\,\sigma(\widetilde x)\big)$.
\smallskip

\begin{Quote}
\begin{lem}\label{lem-asym}
This $~\widetilde{\nu}~$ is an asymptotic measure for the extended
map $~\widetilde{F}:\Sigma\,\times I\,\to\Sigma\,\times\, I$.
\end{lem}
\end{Quote}
\smallskip

{\bf Proof.} We know that almost every orbit $~(\widetilde x_0,\,y_0)\mapsto
(\widetilde x_1,\,y_1)\mapsto\cdots~$ under $~\widetilde{F}~$ converges (in the
sense of Lemma \ref{lem-rto0}) towards the graph of $~\sigma$.
If $~\chi:\Sigma\times I\to\R~$ is any continuous test function, then it 
follows easily that the difference between the time averages
$$\Big( \sum_{0}^{n-1}\chi(\widetilde x_i, \,y_i)\Big)/n\qquad{\rm and}\qquad
\Big(\sum_{0}^{n-1}\chi\big(\widetilde x_i, \,\sigma(\widetilde x_i)\big)
\Big)/n~=~
\Big(\sum_{0}^{n-1}\chi\big(\widehat\sigma(\widetilde x_i)\big)\Big)/n $$
converges to zero as $~n\to\infty$. But the Birkhoff Ergodic Theorem,
applied to the bounded measurable function $~\chi\circ\widehat\sigma:\Sigma
\to\R$, ~ asserts that this last time average converges towards the space 
average
$$  \int\limits_{\Sigma} \chi\circ\widehat\sigma(\widetilde x)\,d\mu_{\Sigma}
(\widetilde x)
~=~\int\limits_{\Sigma\times I}\chi(\widetilde x,\,y)\,d\widetilde{\nu}
(\widetilde x,\,y)\,, $$ 
as required. \qed
\bigskip

{\bf Proof of Theorem \ref{theo-meas} (conclusion).\/} If 
$\,{\mathcal S} f_x > 0\,$ it  follows from Lemma B\ref{lem-prod} and  
Lemma A\ref{A.1.}  that  at least one 
 of the transverse Lyapunov exponents of the boundary circles  is 
strictly positive, hence its corresponding basin has measure zero. If both,
 $\Ly(\A_0)$ and $\Ly(\A_1)$ are  strictly positive, the same lemma implies
that the set of orbits converging to $\A_0 \cup \A_1$ has measure zero,
so that the Lemmas \ref{lem-gam} through \ref{lem-asym} apply.
\medskip

Pushing forward the canonical measure $\mu_{\Sigma}$ on ${\Sigma}$
by the section $\widehat\sigma:
\widetilde{x} \mapsto (\widetilde{x}, \sigma(\widetilde{x}))$, we
obtain an asymptotic measure $\widetilde{\nu}=
\widehat\sigma_*(\mu_{\Sigma})$ for the map $\widetilde{F}$. 
Now, pushing forward again under the projection 
$$ (\widetilde x,\;y)~\mapsto~(x_0,\;y) $$
from $\Sigma\times I$ to $(\R/\Z)\times I=\cC$,
we obtain an $F$-invariant measure $\nu$ on $\cC$.
Since almost every orbit under $\widetilde{F}$ is uniformly distributed with
 respect to $\widetilde{\nu}$, it follows that almost every orbit under $F$ is
 uniformly distributed with respect to $\nu$.
\qed \smallskip

\begin{rem}
In the spirit of  Remark \ref{rem-sepm} one could say, that the 
{\bit separating  measure} for $\widetilde F^{-1}$  is an {\bit asymptotic 
measure} for $\widetilde F$.
\end{rem}

\setcounter{lem}{0}

  \section{Zero Schwarzian}\label{ZS}
This section will study the intermediate case where each diffeomorphism
$~f_x:I\to I~$ has Schwarzian $~{\mathcal S}f_x~$ identically zero. Such a map
is necessarily fractional linear, and can be written for example as
\begin{equation}\label{eq:fl}
  y~\mapsto~\frac{ay}{1+(a-1)y}\qquad{\rm with}\qquad a~>~0\,,
\end{equation}
where $~a~$ is the derivative at $~y=0$.~
It will be convenient to replace $~y~$ by the {\bit Poincar\'e
arclength coordinate}
\begin{equation}\label{eq:pc}
 t(y)~=~\log\rho(0,\,1/2,\, y,\,1)~=~\log\frac{y}{1-y}\,, 
\end{equation}
which varies over the entire real line
for $~0<y<1$,~ with inverse $~y=e^t/(1+e^t)$.~
If we embed the unit interval in the complex open disk of radius $1/2$
centered at $1/2$,~then $~|t|~$ can be described as the distance from the
midpoint, using the Poincar\'e metric for this disk. (Compare Appendix~B.)

Since we are assuming $~{\mathcal S}f_x~$ identically zero, it follows that
each $~f_x~$ preserves cross-ratios or Poincar\'e distances. (See
Equations~(\ref{eq:pc}) and (\ref{eq:pcd}) of Appendix~B.)
Therefore, in terms of the Poincar\'e arclength coordinate $~t$,~ the map
$~f_x~$ will simply be a translation, $~t\mapsto t+c~$ where $~c~$ is a
constant depending on $~x$.~ In other words,
$$ t(f_x(y))=t(y)\;+\;c\,,\qquad{\rm where}\qquad c~=~\log(a)~\in~\R
\qquad{\rm or}\qquad a~=~e^c\,.$$
Using this displacement $~c~$ in place of the
original parameter $~a$~, the 1-parameter group
of fractional linear transformations of the unit interval takes the form
$$  g_c(y)~=~\frac{e^cy}{1+(e^c-1)y}\,,$$ 
where $~g_{c+c'}=g_c\circ g_{c'}$. Given any smooth function $~p~$
from $~\R/\Z~$ to $~\R$,~ we can set $~c=p(x)~$ to obtain
an associated cylinder map
$$   F(x,\,y)~=~ \big(kx,\,g_{p(x)}(y)\big)\,. $$
(The map $~F~$ has an absolutely continuous invariant measure $~dx\,dt$.~
However, this is not very useful since the total area $~\int\!\int dx\,dt~$
is infinite.)\medskip

Using the coordinate $~t\in\R~$ in place of $~y\in(0,1)$,~ the cylinder
map $~F~$ will correspond to the map
$$ (x,\,t)~\mapsto~\big(k\,x,\;t+p(x)\big)$$
of $(\R/\Z)\times\R$.
The dynamics of $~F~$ under iteration is governed by the average
\begin{equation}\label{eq:av}
   A~=~ \int\limits_{\R/\Z} p(x)\,dx
\end{equation}
of the displacement $~p(x)$.
For almost any orbit $~~(x_0,\,t_0)\mapsto\cdots\mapsto (x_n,\,t_n)
\mapsto\cdots$,~~ it follows from the Birkhoff Ergodic Theorem that
the time average
$$ (t_n-t_0)/n~=~\big(p(x_0)+\cdots+p(x_{n-1})\big)/n $$
converges to the space average $~A~$ as $~n\to\infty$. Thus if
$~A>0~$ then it follows that
$~t_n~$ will converge to $~+\infty$.~ In other words,
the corresponding orbit for $~F~$ will converge towards the upper
cylinder boundary $~\A_1$, so that $~\A_1~$ will be a global attractor
under $~F$.~ Similarly, if $~A<0~$ then the lower boundary
$~\A_0~$ will be a global attractor.\medskip

The borderline case where the average 
(\ref{eq:av}) is exactly zero, is much more interesting.
We conjecture that the long term behavior of the sequence of numbers
$t_0,\,t_1,\,t_2,\,\ldots$ is very much like that for a random walk, in which
the successive differences $~\Delta t_n= t_{n+1}-t_n~$ are identically
distributed independent random numbers with mean zero.
 In particular, we believe that the
following theorem will be true whenever the periodic function $~p(x)~$ has
average zero. However, the proof will apply only in the following
very special case.

\begin{quote} {\bf Hypothesis.} We now assume that $~p(x)~$ is a step function
which takes a constant value on each of the $~k~$ intervals
 $~j/k\le x<(j+1)/k.$
\end{quote}

\noindent If $x_0$ is randomly chosen, it
 then follows easily that the successive steps
$$ \Delta t_n~=~t_{n+1}-t_n~=~p(x_n) $$
actually are identically distributed independent random variables, which
 do not depend on the value of $~t_n$.~ (In fact
$~\Delta t_n~$ depends only on the $~n$-th entry in the base $k$ expansion
of $~x_0$.) Choose some
number $~N\ge 0~$ and define three sequences as follows. Let
\smallskip

\hskip 1cm $a_n~$ be the number of integers $~1\le i\le n~$ with\hskip .75cm
 $t_i>N$,\smallskip

\hskip 1cm $b_n~$ the number of such integers with\hskip 2.25cm
 $~|t_i|\,\le\,N$, and\smallskip

\hskip 1cm $c_n~$ the number of such integers
with\hskip 2.3cm $~t_i\,<\,-N$.\medskip

\noindent Thus the associated frequencies  $~a_n/n,~b_n/n,~ c_n/n~$ will have
 sum equal to $~+1$.

\begin{Quote}
\begin{theo}\label{th:Schw=0}
Let $~p(x)~$ be a step function as described above, not identically zero
but with average $~\int p(x)\,dx~$ equal to zero.
Then for arbitrary $t_0$ and for Lebesgue almost every
 $x_0$, the ratios $~b_n/n~$
associated with the orbit $~(x_0,\,t_0)\mapsto (x_1,\,t_1)\mapsto\cdots$~
will converge to zero as $~n\to\infty$, but $~a_n/n~$ and
$~c_n/n~$ will not tend to any limit. In fact, these ratios
 vary so wildly that\vspace{-.2cm}
$$ \liminf\,(a_n/n)~=~0~<~\limsup\,(a_n/n)~=~1\,,$$
with a similar statement for $~c_n/n$.
\end{theo}
\end{Quote}
\smallskip

 In terms of the original cylinder map, this means that most orbits
spend most of the time extremely close to one or the
other of the two boundaries; very occasionally jumping
 from one boundary to the other but
in such an irregular way that there is no limiting asymptotic measure.
\medskip

 We are indebted to Harry Kesten and especially to
Mikhail Lyubich for very substantial
 help with the proof. To begin the argument,
note that for each $t_0$ and each $~x_0\in{\mathbb R}/{\mathbb Z}~$ the number
\break $~L=\limsup_{n\to\infty}( a_n/n)~$ is well defined, and is invariant
under finite permutations of the sequence of differences $~\Delta t_n$, or
equivalently under finite permutations of the entries in the base $k$ expansion
of $x_0$. Therefore, according to the Hewitt-Savage {\bit Zero-One Law\/},
each set
$$ \{x_0~~;~~ L~\le~{\rm constant} \}$$
has measure either zero or one.  (See \cite[IV.6]{F2}.)
It follows easily that $~L~$ takes some constant value for almost all $~x_0$.

 First consider the case $N=0$.
We will prove that $~L=1~$ except on a set of measure zero. Choose some
small $~\epsilon>0~$ and let $X_n(\epsilon)~$ be the set of all $~x_0~$
for which $~a_n/n>1-\epsilon$.~ Evidently the intersection
$$\bigcap_n\Big(X_n(\epsilon)\cup X_{n+1}(\epsilon)\cup X_{n+2}(\epsilon)\cup
\cdots\Big) $$
is precisely the set of $~x_0~$ for which $~L=\limsup(a_n/n)~$
 satisfies  $~L>1-\epsilon$. Thus, if we assume that $~L\le 1-\epsilon~$ on a
set of $~x_0~$ of positive (and hence full) measure, then it would certainly
follow that the measure of
$~X_n(\epsilon)~$ must tend to zero as $~n\to\infty$.~ But according to the
{\bit Arcsine Law\/}, the measure of $~X_n(\epsilon)~$ converges to the limit
$$  1-\frac{2}{\pi}\,\arcsin\sqrt{1-\epsilon}~~=~~
 \frac{2}{\pi}\,\arcsin\sqrt\epsilon~~>~~0 $$
as $~n\to\infty$.~ (See \cite[III.4; and  {\bf 2}, 1966, XII.8]{F1}.)
 This contradiction proves that $~\limsup(a_n/n)=1~$
almost everywhere. Since $~a_n+c_n\le n$,~ it follows that 
$~\liminf\,(c_n/n)=0$;
and a similar argument with $~a_n~$ and $~c_n~$ interchanged completes
 the proof for the case $~N=0$.

To complete the proof of Theorem \ref{th:Schw=0},
 choosing some arbitrarily large $~N$, we must show
that the associated sequence $~b_n/n~$ converges to zero for almost all $~x_0$.
The proof will be based on the following.\medskip

Consider a random walk $~\tau_0,\,\tau_1,\,\ldots~$ on the circle ${\mathbb R}/
{\mathbb Z}$,~ starting with some specified $~\tau_0$,~
where the differences $~\Delta\tau_n=\tau_{n+1}-\tau_n~$ are
identically distributed independent random variables which do not depend
on $~\tau_n$.
\smallskip

\begin{Quote}\begin{lem}\label{eq-distr} The resulting sequence
$~\{\tau_n\}~$ is uniformly distributed around the circle with probability $+1$
if and only if the following condition is satisfied:

{\bf (*)} There is no finite
cyclic subgroup $~G\subset{\mathbb R}/{\mathbb Z}~$ such that
 $~\Delta\tau_n\in G~$ with probability $+1$.
\end{lem}\end{Quote}
\smallskip

{\bf Proof.} (Compare \cite{Levy}, \cite{KI}.)
Clearly this subgroup condition is necessary. To prove that
it is sufficient, let $~\mu~$ be the common probability distribution for the
differences $~\Delta\tau_n$,~ let $~S\subset\R/\Z~$ be the support of $~\mu$,~
and let $~S^{\mathbb N}=S\times S\times S\times\cdots~$ be the space of
sequences ${\bf s}=(s_0,\,s_1,\,\ldots)~$ of elements of $~S$,~ provided with the
shift invariant measure $~\mu^{\mathbb N}=
\mu\times\mu\times\cdots$. It is not hard to see that the skew product map
$$ F\big(\tau,\,s_0,\,s_1,\,\ldots)~=~(\tau+s_0,\,s_1,\,s_2,\,\ldots) $$
from $~(\R/\Z)\times S^{\mathbb N}~$ to itself preserves the measure
$~\lambda\times\mu^{\mathbb N}$,~ where $~\lambda~$ is Lebesgue measure.
We will prove that $~F~$ is ergodic; or equivalently that:

\begin{quote} \it Every bounded
measurable $~F$-invariant function  $~\phi:(\R/\Z)\times S^{\mathbb N}~\to
~\R~$ is constant almost everywhere.\end{quote}

\noindent Given such a function $~\phi~$, let
$$ \phi_n(\tau,\,s_0,\,\ldots,\,s_{n-1}) ~=~\int\phi(\tau,\, {\bf s})
\,d\mu^{\mathbb N}(s_n,\,s_{n+1},\,\ldots)$$
be the average of $~\phi(\tau,\, {\bf s})~$
 over all possible choices of $~s_n,\,s_{n+1},\,\ldots\,$.~ Clearly
 $~\phi_n~$ can also be expressed as an average over all choices of $~s_n$,
$$ \phi_n(\tau,\,s_0,\,\ldots,\,s_{n-1}) ~=~\int\phi_{n+1}(\tau, s_0,\,\ldots,
\,s_n)\,d\mu(s_n)\,.$$
On the other hand, using the condition $~\phi=\phi\circ F~$ of $~F$-invariance,
we see that $~\phi_n(\tau,\,s_0,\,\ldots,\,s_{n-1})~$ can also be described as
the average  of $~\phi(\tau+s_0,\,s_1,\,s_2,\,\ldots)$ over all choices
of $~(s_n,\,s_{n+1},\,\ldots)~$.~ But by definition, this is equal to
 $~\phi_{n-1}(\tau+s_0,\,s_1,\,\ldots,\,s_{n-1})$.
 Thus, inductively, it follows that
\begin{eqnarray}\label{comp-phi-n}
\phi_n(\tau,\,s_0, \ldots,\,s_{n-1})&=&\phi_{n-1}(\tau+s_0,\,s_1,\,\ldots,
\,s_{n-1})\,=\,
\phi_{n-2}(\tau+s_0+s_1,\,s_2,\,\ldots,\, s_{n-1})\,=\,\cdots~~~\nonumber \\
&=&\phi_0(\tau+s_0+\cdots+s_{n-1})\,.
\end{eqnarray}
In particular, note that
\begin{equation}\label{conv}
\phi_0(\tau)~=~\int\phi_1(\tau,\,s)\,d\mu(s)~=~\int\phi_0(\tau+s)\,d\mu(s)\,.
\end{equation}

If the condition {\bf(*)} is satisfied, then we will use this last equation
to prove that $~\phi_0~$ is constant almost everywhere. Let
$$ \widehat\mu(q)~=~\int {\bf e}(-qs)\,d\mu(s) $$
be the Fourier transform of $~\mu$,~ where $~q\in\Z$,~ and where
$~{\bf e}(t)~$ is an abbreviation
for $~e^{2\pi it}$. Thus $~\widehat \mu(q)~$ is a weighted average of points
on the unit circle; hence $~|\widehat\mu(q)|\le 1$.~ For $~q\ne 0$,~ the
condition {\bf(*)} guarantees that the weight is not all concentrated at
points $~s~$ such that $~qs\equiv 0~({\rm mod}~\Z)$,~
so it follows that $~\widehat\mu(q)\ne 1$.\smallskip

The Fourier transform $~\widehat\phi_0(q)~$ is defined similarly as the
 integral of $~{\bf e}(-q\tau)\phi_0(\tau)\,d\tau$.~ Recall that a bounded
measurable function on the circle is constant almost everywhere if and
only if its $~q$-th Fourier coefficient is zero for every $~q\ne 0$.
 Multiplying equation 
(\ref{conv}) by $~{\bf e}\,(-q\tau)\,d\tau~$ and then integrating, using the
 substitution $~\eta=\tau+s$,~ we obtain
$$ \widehat\phi_0(q)~=~\int{\bf e}(-q\eta)\,\phi_0(\eta)\,d\eta
\,\int{\bf e}(qs)\,d\mu(s)~=~\widehat\phi_0(q)\,\widehat\mu(-q)\,.$$
For every $~q\ne 0$,~ since $~\widehat\mu(-q)\ne 1$,~
it follows that $~\widehat\phi_0(q)=0$.~ Therefore $~\phi_0~$ takes some
constant value $~v~$ almost everywhere. Using equation (\ref{comp-phi-n}),
it follows that every $~\phi_n~$ also takes the value $~v~$ almost everywhere.

We can now prove that $\phi~$ is constant almost everywhere. Suppose to the
contrary, for example, that $~\phi(\tau,\,{\bf s})>v+\epsilon>0~$ on a set
$~\Sigma~$ of positive measure. Choose a point of density
 $~(\tau^*,\,{\bf s}^*)~$ for $~\Sigma~$. Since $~\phi~$ is bounded,
 it would follow that the average of $~\phi~$ over a small neighborhood of
 $~(\tau^*,\,{\bf s}^*)~$ is strictly greater than $~v$.~ But the fact that
 each $~\phi_n~$ equals $~v~$ almost everywhere implies that every such
 average is also equal to $~v$.~
 Thus $~\phi~$ must be constant almost everywhere; which proves ergodicity.

Now the Birkhoff Ergodic Theorem implies that almost every orbit of $~F~$
is uniformly distributed with respect to the measure
 $~\lambda\times\mu^{\mathbb N}$.~ Projecting to the first coordinate, it
follows that almost every sequence
$$ \tau,~~\tau+s_0,~~\tau+s_0+s_1,~~\ldots $$
is uniformly distributed with respect to $~\lambda$. This
 completes the proof of Lemma \ref{eq-distr}.\qed\bigskip

{\bf Proof of Theorem \ref{th:Schw=0} (conclusion).} Recall that $~b_n~$ is
the number of $0\le j<n~$ with $~|t_j|\le N\,.$ 
We must show that the ratio $~b_n/n~$ tends to zero
with probability $+1$ as $~n\to\infty\,.$ Choosing some large number
$~L\gg N$,~ the quotients
$~t_n/L~~ ({\rm mod}~{\mathbb Z})~$ form a random walk on the
circle. If $~L~$ is chosen so that
 the ratio $~\Delta t_n/L$ is irrational with probability $\,>0$,~
then by Lemma \ref{eq-distr}, these quotients are uniformly distributed around
the circle. Let $~b_n(L)\ge b_n~$ be the number of $~j\in[0,\,n)~$
with $~t_j~$ congruent to
an element of $~[-N,\,N]~$ modulo $L{\mathbb Z}$.~ Then
 it follows that $~b_n(L)/n~$
converges to $~2N/L$ as $~n\to\infty~$. Since $~L~$ can be arbitrarily large,
 this proves that $~b_n/n\to 0$, as required.\qed

\appendix

\renewcommand{\theequation}{A\arabic{equation}} 

\setcounter{equation}{0}
\setcounter{lem}{0}
\section*{Appendix A: The Transverse Exponent.}

Let $(x,\,\iota)$ be any point of the boundary circle $\A_\iota$, where $\iota$
 can be either $0$ or $1$.
By definition the {\bit transverse Lyapunov exponent} along the
circle $\A_\iota$ at $(x,\,\iota)$ is defined by
$${\rm Lyap}_{{\mathcal A}_{\iota}}(x)=\lim_{k \rightarrow \infty}\frac{1}{k}
\log \Big\vert
\frac{ \partial F^{\circ k}}{\partial y}(x,y)\Big\vert\qquad
{\rm evaluated~at}\quad y=\iota\,,
$$
whenever this limit exists.
(In this case, {\bit transverse} really means {\bit normal}.)
Here $~F(x,y)=\big(kx,~f_x(y)\big)~$ as usual.
By the chain rule, the above expression can be written as
$${\rm Lyap}_{{\mathcal A}_{\iota}}(x)=\lim_{k\rightarrow \infty}\frac{1}{k}
\log\Big( f'_{x_0}(\iota)f'_{x_1}(\iota)\cdots f'_{x_{k-1}}(\iota)\Big)$$
where $x_0 \mapsto x_1\mapsto \cdots$ is the orbit of  $x = x_0$.\smallskip

Let us denote by $\lambda$ the 2-dimensional Lebesgue measure on the cylinder
 ${\mathcal C}=(\R/\Z)\times I$, and by $\lambda_x$ the 1-dimensional Lebesgue
 measure along $\R/\Z$.  Since  $\lambda_x$ is ergodic and invariant under 
multiplication by $k$ (see Equation~(\ref{eq-for})), it follows from the
 Birkhoff Ergodic Theorem that  this transverse Lyapunov exponent is
defined and independent of $x$ for almost all
$x$, and is equal to the integral
$${\rm Lyap}_{{\mathcal A}_{\iota}}= 
\int\limits_{\R/\Z}\log \Big(f'_{x}(0)\Big) dx,$$
for almost all $x$.
\smallskip

Let us prove now the Lemma stated in  \S \ref{Preli}. 
\smallskip

\begin{Quote}
\begin{alem}\label{A.1.}
{\it For $\iota$ equal to zero or one, let $\B_{\iota}$ be the
attracting basin of the circle $\A_\iota$. If the transverse Lyapunov exponent
\begin{equation}
\Ly(\A_{\iota})=\int\limits_{\R/\Z} \log\left(f'_x(\iota)\right)dx
\end{equation}
is negative, then the basin $\B_{\iota}$ has strictly positive measure.
In fact, for almost every $x\in \R/\Z$ the basin $\B_\iota$ intersects
the ``fiber'' $x\times I$ in an interval of positive length. On the other
hand, if $~\Ly(\A_{\iota})>0~$ then $\B_{\iota}$ has measure zero.}
\end{alem}
\end{Quote}
\smallskip

{\bf Proof.\/} First consider the case $\iota=0$ with $\Ly(\A_0)<0$.
By Taylor's expansion restricted to the fiber over $x$, we have
$$f_x(y)= f_x'(0)y + O(y^2)\, ,$$
uniformly for all $\,(x,y)\, \in\, (\R/\Z)\,\times\, I$.  Choose $K>0$ so that,
$$f_x(y)~\leq~ y \big(f_x'(0)+Ky\big)\qquad {\rm for\quad all}\qquad (x,y)\,.$$
For any $\eta>0$, it follows that
\begin{equation}\label{eq-ubf}
f_x(y)~\leq~ y\big(f_x'(0)+\eta\big)
 \qquad {\rm whenever} \qquad y < \frac{\eta}{K}\,.
\end{equation}
Since $\Ly({\mathcal A}_0)<0$, we can choose $\eta>0$ small enough so that
\begin{equation}\label{eq-ni}
\int\limits_{\R/\Z}\log\big(f'_x(0)+\eta\big)dx~<~ 0\,.
\end{equation}
It will be convenient to introduce the notation
\begin{equation}\label{eq-axj}
a(x)~=~ \log\big(f_{x}'(0)+\eta\big)\,.
\end{equation} 
Consider some orbit  $~(x_0,\,y_0) \mapsto(x_1,\, y_1) \mapsto(x_2,\, y_2)
 \mapsto \cdots$.
By the Birkhoff Ergodic Theorem, the averages
$$\frac{1}{n}\Big(a(x_0)+a(x_1)+\ldots+ a(x_{n-1})\Big) $$
converge to $~\int_{\R/\Z}\, a(x)\,dx<0~$ for almost all $x_0$. In particular,
 it follows that the $n$-fold sum
$$  A_n(x_0) ~=~a(x_0)+a(x_1)+\ldots+ a(x_{n-1}) $$
converges to negative infinity as $n\to\infty$. Hence the maximum
$$A_{\rm max}(x)~=~\max_{n\ge 0} A_n(x) $$ 
is certainly defined and finite for almost all $x$.
 Now suppose that
\begin{equation}\label{up-bd}
 y_0~\le~ \frac{\eta}{K}\,e^{-A_{\rm max}(x_0)}\,.
\end{equation}
Then a straightforward induction shows that
$$ y_n~\le~
 \frac{\eta}{K}e^{A_n(x_0)-A_{\rm max}(x_0)}~\le
~ \frac{\eta}{K} $$
for all $n$. Since $A_n(x_0)$ converges to $-\infty$,
 it follows that $y_n$ tends to zero, so
that $(x_0,\,y_0)$ belongs to the attracting basin $\B_0$. Since the
right side of the inequality (\ref{up-bd}) is a measurable function of $x_0$,
defined and strictly positive almost everywhere, it follows that its integral
is strictly positive. Evidently this integral is a lower bound for the area
of $\B_0$. Thus $\B_0$ has positive measure as required.\medskip

The proof for the case $\Ly(\A_0)>0$ is completely analogous. However, it
requires us to make use of the hypothesis that $f_x'(y)$ is strictly positive,
even for $y=0$, so that we can choose a small $\eta$ with $0<\eta<f_x'(0)$
everywhere, and with
\begin{equation}\label{posint}
\int\log\big(f'_x(0)-\eta\big)\,dx~>0\,.
\end{equation}
The estimate (\ref{eq-ubf}) is then replaced by
\begin{equation}\label{lb}
  f_x(y)~\ge~ y\big(f'_x(0)-\eta\big)\qquad{\rm whenever}\qquad y
<\frac{\eta}{K} \,.
\end{equation}
 Now suppose that the basin $\B_0$ has positive measure. Then, for a set of
$x_0$ of positive measure, we could find orbits
 $(x_0,\,y_0)\mapsto(x_1,\,y_1)\mapsto \cdots$ which satisfied
 $0<y_n<\frac{\eta}{K}$
for all $n$. But using (\ref{posint}) and
(\ref{lb}) it is not hard to see that this is impossible.
Therefore $\B_0$ has measure zero. The arguments for the basin $\B_1$
are completely analogous. \qed

\renewcommand{\theequation}{B\arabic{equation}} 

\setcounter{equation}{0}
\setcounter{lem}{0}
\section*{Appendix B:  Schwarzian Derivative and Cross-Ratios.}
Recall that the Schwarzian derivative $\,{\mathcal S}f\,$ of a $~C^3~$
interval
diffeomorphism $~f: I \rightarrow I~$ was defined in Equation~(\ref{eq-ns}).
The statement $~{\mathcal S}f < 0~$ will mean that the inequality ${\mathcal S}
 f(y) < 0$ holds for $~y~$ in a  dense open subset of  $I$; and similarly
for  $\,{\mathcal S} f > 0\,$ (or $\,{\mathcal S} f = 0\,$).
\smallskip

\begin{Quote}
\begin{bprop}\label{prop-propt}
{\it The Schwarzian derivative has the following properties:}

{\bf 1.}~~{\it The sign of $\,{\mathcal S}f\,$ is preserved under iteration of
 $f$. For example if $\,{\mathcal S}f\,<\,0\,$ and $\,{\mathcal S}g\,<0\,$,
 then $\,{\mathcal S}(f\circ g)\, < \,0\,$.}

{\bf 2.}~~{\it $\,{\mathcal S}f\,<0\,$ if and only if $\,{\mathcal S}f^{-1}\,
 >0\,$.}

{\bf 3.}~~{\it $\,{\mathcal S}f\, <\,0\,$ if and only if the function 
$~\varphi(y)=1/\sqrt{\vert f'(y) \vert}~$ is concave  upwards.}

{\bf 4.}~~{\it $\,{\mathcal S}f\, =\,0\,$ if and only if  $f$ is a fractional 
linear transformation, $f(x)=(ax+b)/(cx+d)$ where we may assume that  
$ad-bc~=~\pm~1$.}
\end{bprop}
\end{Quote}
\smallskip

{\bf Proof.\/}  

\begin{enumerate}
\item A straightforward calculation shows that the Schwarzian derivative
of a composition is given by the formula
\begin{equation}\label{sch-comp}
{\mathcal S}(f\circ g)~=~ (g')^2\,{\mathcal S}f\, +\, {\mathcal S}g\,,
\end{equation}
and the conclusion follows easily.
\item This follows by taking $g=f^{-1}$ in equation (\ref{sch-comp}) and noting
that the identity map has Schwarzian zero.
\item It is not hard to calculate that the second derivative of $~\varphi(y)=
 1/\sqrt{\vert f'(y) \vert}~$ satisfies the equation
\begin{equation}\label{phi''}
{\mathcal S}f~ = ~-2\varphi''(x)/\varphi(x)\,.
\end{equation}
Thus ${\mathcal S}f(y)<0$ on a dense open set if and only if $\varphi''(y)>0$
on a dense open set, and the assertion follows.
\item From equation (\ref{phi''}) we see that the Schwarzian is zero if
and only if the function $\,\varphi(y)\,$ is linear, say $\;\varphi(y)
=cy+d\;$ or in other words
$$ f'(y)~=~\pm 1/(cy+d)^2\,.$$
Integrating, we see that this is true if and only if  $f$ is fractional linear.
\qed
\end{enumerate} 
\smallskip

Now consider a $C^3$-diffeomorphism $f:I\to I'$ where $I'$ may be a different
interval of real numbers.
 A fixed point $~y=f(y)\in I\cap I'~$ will be called {\bit strictly
attracting\/} if $~|f'(y)|<1~$ and {\bit strictly repelling\/} if $~|f'(y)|>1$.
\smallskip

\begin{Quote}
\begin{blem}\label{3fp}
{\it Suppose that ${\mathcal S} f<0$ throughout a dense open subset of $I$.
If $f$ has two fixed points, then at least one must be strictly attracting.
 Furthermore, if there are three fixed points
then the middle one must be strictly repelling and the other two must be
strictly attracting. Such a map can never have four fixed points.}
\end{blem}
\end{Quote}
\smallskip

{\bf Proof.} Given two fixed points $\alpha<\beta\,$,~ we will first show that
at least one of the two must be strictly attracting. In fact, it follows from
the Mean Value Theorem that some point $\alpha<x<\beta$ must satisfy $f'(y)=1$.
If both $f'(\alpha)\ge 1$ and $f'(\beta)\ge 1$, then the graph of $f'$ would
have to have
have a local minimum somewhere in the open interval $(\alpha,\,\beta)$.
Hence the function $\varphi(y)$ would have a local maximum. Since
$\varphi''(y)>0$ on a dense open set, this clearly leads to a contradiction.

Now consider three fixed points $\alpha<\beta<\gamma$. Then The Mean Value
Theorem yields points $y\in (\alpha,\beta)$ and $y'\in(\beta,\gamma)$ with
$f'(y)=f'(y')=1$. Again, if $f'(\beta)\le 1$, then the graph of $f'$ would have
a local minimum, yielding a contradiction. Therefore $\beta$ is strongly
repelling, hence $\alpha$ and $\gamma$ must be strongly attracting. Evidently
this leaves no possibility for a fourth fixed point.\qed
\smallskip

\begin{Quote}
\begin{blem}\label{lem-prod}
{\it If $\,{\mathcal S}f\, <\, 0\, $ for an orientation preserving 
diffeomorphism $f:I \rightarrow I$, then   $f'(0)f'(1)< 1$.}
\end{blem}
\end{Quote}
\smallskip

{\bf Proof.\/} By the previous lemma, at least one of the two boundary fixed
points must be strictly attracting: either $~f'(0)<1$, or $~f'(1)<1$, or both.
Now consider the auxiliary function $$g(y)~=~1-f(1-y)\,.$$
 Evidently $g$ also has
negative Schwarzian, with $g'(0)=f'(1)$ and $f'(0)=g'(1)$. The composition
$f\circ g$ also has negative Schwarzian; hence $f\circ g$ has derivative
less than one at at least one of the two endpoints. But the derivative at
either endpoint is equal to the product $f'(0)f'(1)$. It follows that this
product is less than one, as required.\qed\medskip

{\bf Definition.\/} The {\bit cross-ratio} of four distinct real numbers
will mean the expression
\begin{equation}\label{eq-cr}
\rho(y_0, y_1, y_2, y_3)= \frac{(y_2-y_0)(y_3-y_1)}{(y_1-y_0)(y_3-y_2)}.
\end{equation}
(The reader should take care, since conflicting notations are often used.)
Note that
$$\rho(y_0,\,y_1,\,y_2,\,y_3)>1\qquad{\rm whenever}\qquad y_0<y_1<y_2<y_3\,.$$
Evidently the cross ratio remains invariant whenever we replace each $y_i$
by $ay_i+b$ with $a\ne 0$. A brief computation shows that it also remains
invariant when we replace each $y_i$ by $1/y_i$. Since every fractional linear
transformation can be expressed as a composition of affine maps and inversions,
it follows that the cross-ratio is invariant under fractional linear
transformations.
\medskip

We will say that a monotone map $f$ {\bit increases cross-ratios\/} if
$$ \rho\big(f(y_0),\,f(y_1),\,f(y_2),\,f(y_3)\big)~>~
 \rho(y_0,\,y_1,\,y_2,\,y_3)\qquad{\rm whenever}\qquad y_0<y_1<y_2<y_3\,.$$

\begin{Quote}
\begin{blem}\label{lem-nsi}{\bf (Allwright.)\/} 
{\it Again let $f:I\to I'$ be a $C^3$-diffeomorphism.
Then $f$ increases cross-ratios if and only if $~{\mathcal S} f<0~$ throughout
some dense open subset of $I$.}
\end{blem}
\end{Quote}

{\bf Remark.} In  Lemmas B\ref{3fp}, B\ref{lem-prod}, and
 B\ref{lem-nsi}, note that we can obtain
a corresponding statement for the case $~{\mathcal S} f>0~$ simply by applying
the given statement to the inverse map from $I'$ to $I$. For example: {\it
$f$ decreases cross-ratios if and only if $~{\mathcal S} f>0~$ on a dense open
set.}\medskip

{\bf Proof of Lemma B\ref{lem-nsi}.} (Compare \cite{A}.) First suppose
that $~{\mathcal S}f<0~$ on a dense open set.
Given points $~y_0<y_1<y_2<y_3\,$, after composing $f$ with a fractional
linear transformation, we may assume that $~f~$ fixes the three points
$~y_0,\,y_1,\,y_3$. If $~{\mathcal S} f<0~$
on a dense set, then $~y_0,\,y_3~$ are attracting and $~y_1~$ is repelling
by Lemma
B\ref{3fp}. Since there can be no fixed point between $y_1$ and $y_3$, it
follows that $f$ moves every intermediate point to the right. Thus
$f(y_2)>y_2\,$, and it follows easily that $f$ increases the cross-ratio
$\rho(y_0,\,y_1,\,y_2,\,y_3)$.

Conversely, if ${\mathcal S} f$ is not negative on a dense open set, then it
must either be strictly positive somewhere, or identically zero on some
interval. In the first case, it would decrease some cross-ratio, and in the
second case it would be fractional-linear and hence preserve cross-ratios
within this interval. This completes the proof.
\qed\medskip

\begin{brem} If $0 < y_1 <y_2 < 1$, then the {\bit Poincar\'e
distance} between $y_1$ and $y_2$ within $(0,1)$ can be defined as
\begin{eqnarray}\label{eq-ip}
d_{[0,1]}(y_1 , y_2)&=&\int\limits_{y_1}^{y_2}\Big(\frac{1}{y} + 
\frac{1}{(1-y)}
\Big) dy \nonumber\\
&=&\log \Big(\frac{y_2(1-y_1)}{y_1(1-y_2)}\Big)
~=~\log\;\rho(0,y_1,y_2,1)\,. \nonumber
\end{eqnarray}
This can be identified with the usual Poincar\'e distance within a
complex disk having the interval $[0,1]$ as diameter. In terms of the 
{\bit Poincar\'e arclength coordinate} of Equation~(\ref{eq:pc}), the 
Poincar\'e distance formula can also be written as,
\begin{equation}\label{eq:pcd}
d_{[0,1]}(y_1 , y_2)=|t(y_2)-t(y_1)|.
\end{equation}
\end{brem}

\vspace{.6cm}

\parbox{8cm}{
{\sc Araceli Bonifant,\\
Department of Mathematics,\\
University of Rhode Island,\\
Kingston, RI. 02881-0816}\\
{\it E-mail address:} bonifant@math.uri.edu\\}
\parbox{8cm}{
{\sc John Milnor,\\
Institute for Mathematical Sciences,\\
Stony Brook University,\\
Stony Brook, NY. 11794-3660.\\}
{\it E-mail address:} jack@math.sunysb.edu\\}

\end{document}